# Research on FAST Active Reflector Adjustment Algorithm Based on Computer Simulation

Tongyue Shi[1], Siyu Tao[1], Haining Wang[2]

(1 School of Computer Science and Technology, Soochow University, Suzhou, China)
(2 School of Mathematical Sciences, Soochow University, Suzhou, China)
`{tyshi, sytao, hnwang}@stu.suda.edu.cn`

## Abstract

Based on the background of the 2021 Higher Education Club Cup National College Students Mathematical Contest in Modeling A, according to the relevant data of the China Sky Eye (FAST) radio telescope, the main cable nodes and actuators are adjusted and controlled by mathematical modeling and computer simulation methods to realize the active reflector. The adjustment of the shape enables it to better receive the signal of the external celestial body and improve the utilization rate of the reflector, so as to achieve a higher receiving ratio of the feed cabin. In this paper, a point set mapping algorithm based on the rotation matrix of the spatial coordinate axis is proposed, that is, the mapping matrix is obtained by mathematical derivation, and the linear interpolation algorithm based on the original spherical surface and the ideal paraboloid to solve the working paraboloid is obtained. When the interpolation ratio is adjusted to 89% to satisfy the optimal solution under realistic constraints. Then, a three-dimensional spatial signal reflection model based on spatial linear invariance is proposed. Each reflective panel is used as an evaluation index. For each reflective surface, the 0-1 variable of signal accessibility on the feeder is defined. The signal is mapped to the plane where the feeder is located, which reduces a lot of computational difficulty. Finally, it is found that the panel of the feed cabin that can receive the reflected signal accounts for 19.3%. Compared with the original spherical reflection model, the working paraboloid model established in this paper has a The signal ratio has increased by 224%.

Keywords: mathematical modeling; computer simulation; point set mapping and surface interpolation; spatial linearity invariance; three-dimensional reflection law

# 1 Research Background

The Five-hundred-meter Aperture Spherical Radio Telescope (FAST) ( Figure 1) is a major infrastructure construction project in my country's National "Eleventh Five-Year Plan". The completion of FAST has made our country in the forefront of the world in the field of astronomical detection, and maintains a leading position [1].

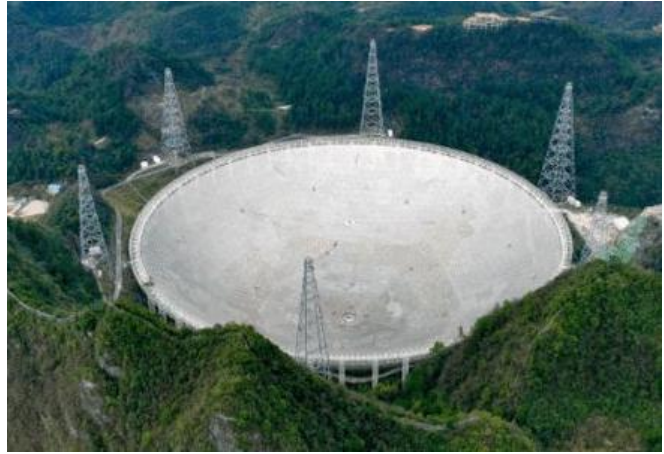

Figure 1 Aerial view of the China Sky Eye FAST radio telescope

The main structure of FAST includes active reflector, signal receiving system and other related auxiliary, load-bearing, measurement and other systems (see Figure 2). The precise structure and observation capability of the telescope. Here, it is particularly important to study the mutual cooperation of the various components of FAST. Usually, it is necessary to determine a parabolic deformation scheme that is conducive to the reflection and collection of signals by the radio telescope according to the azimuth angle of the celestial body. The specific method is to determine an ideal paraboloid to meet the adjustment constraints of the reflector plate of the FAST radio telescope, and then adjust the radial expansion and contraction of the actuator to make the reflector a paraboloid in the working state and fit the paraboloid in the ideal state as much as possible. In this way, a better reflection and collection effect of celestial electromagnetic wave signals can be obtained.

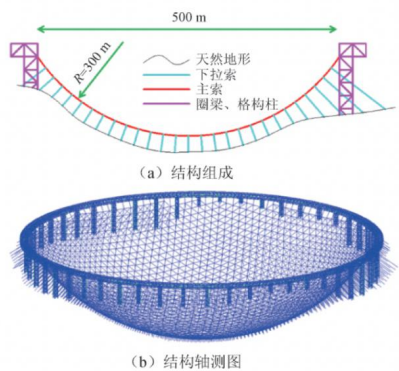

Figure2 The overall structure of the China Sky Eye FAST radio telescope[2]

## 1.1 Questions raised

There are various structures under the main structure of the FAST radio telescope. Among them, the active reflection system composed of the main components such as the main cable net, the reflection panel, the pull-down cable, the actuator and the supporting structure is an adjustable spherical surface (see Figure 3). It can be divided into two states: the reference state and the working state. The task of adjusting the reflective surface into a working paraboloid is mainly accomplished by adjusting the pull-down cable to cooperate with the actuator.

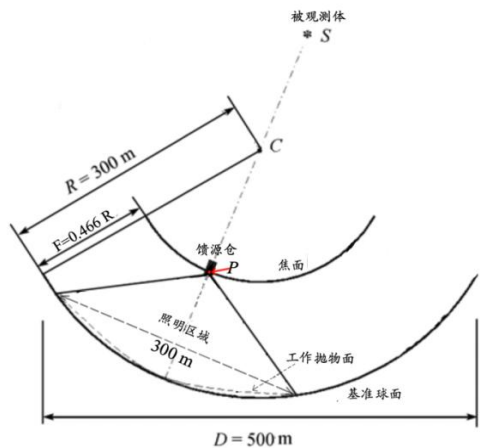

Figure 3 Schematic diagram of the cross-section of the China Sky Eye FAST radio telescope

According to relevant data and information, the problems related to radio telescopes are summarized into the following three problems:

1) When the celestial object to be measured is located directly above the reference sphere, an ideal paraboloid needs to be determined by taking into account the adjustment factors of the reflective panel.

2) When the celestial object to be measured is located at a certain angle, an ideal paraboloid needs to be determined, and an adjustment model of the reflection panel is established to adjust the expansion and contraction of the relevant actuators, so that the reflection surface is as close as possible to the obtained ideal paraboloid. And find the vertex coordinates of the ideal paraboloid, and adjust the results of the main cable node number, position coordinates, and the expansion and contraction amount of each actuator within the 300-meter diameter of the rear reflection surface.

3) Based on the reflection surface adjustment scheme in (2), it is necessary to calculate and adjust the ratio of the reflected signal received in the effective area of the feed cabin to the reflection signal of the reflection surface within the 300-meter aperture (the receiving ratio of the rear feed cabin), and then It is compared with the receiving ratio of the reference reflective sphere, and finally the corresponding conclusions are drawn.

## 2 Model Assumptions

Assume that in the reference state, all main cable nodes are located on the reference sphere.

It is assumed that all structures are normally stressed during the change process, and no large deformation will occur and the overall structure will be damaged.

Assume that each reflective panel is part of the reference sphere and is non-porous.

Assume that both the electromagnetic wave signal and the reflected signal travel in a straight line.

Assume that the radial expansion and contraction of the top end of the actuator is 0 in the reference state.

Assume that the radial retractable range of the actuator is -0.6 to +0.6 meters.

All radial movement directions are assumed to be along the direction of the actuator, ie along the radius of the reference sphere.

Neglect the impact of natural disasters such as landslides, earthquakes, rainstorms and other factors on the main structure of FAST.

Assume that the signal of the celestial body observed by FAST is parallel to the line connecting the observed celestial body and the center of the reference sphere.

## 3 Symbol description

Table 1 Symbol description comparison table

| symbol | explanation | unit |
|---|---|---|
| $id$ | The number of the main index node | / |
| $M_{x,y,z}$ | X, Y, Z coordinates of the main index node | / |
| $D_{x,y,z}$ | Actuator lower endpoint (ground anchor point) X, Y, Z coordinates | / |
| $U_{x,y,z}$ | The X, Y, Z coordinates of the upper end point (top) in the reference state | / |
| $a, c$ | Ideal Paraboloid Quadratic Coefficients and Constants | / |

| $N_{x,y,z}$ | The corresponding point when the main cable node moves to the ideal paraboloid | / |
| $R_{x,y,z}$ | The point to which the main cable node actually moves | / |
| $\boldsymbol{n}_z$ | Incident signal direction vector | / |
| $\boldsymbol{n}_0$ | Reflector normal vector | / |
| $\boldsymbol{n}_{x,y,z}$ | Reflected Signal Direction Vector | / |
| $R$ | Radius of the base sphere | m |
| $F$ | The difference between the radius of the focal plane and the radius of the two concentric spheres of the reference sphere | m |

# 4 PROBLEM ANALYSIS

## 4.1 Analysis of Problem 1

The condition of problem 1 is that the celestial object S to be measured is located directly above the reference sphere. In this case α = 0°, β = 90°, it is easier to determine the ideal paraboloid by considering the reflective panel factor. This is an ideal state. At this time, the three points of the observed celestial body S, the reference sphere center C and the feed cabin P are on the same straight line. The three-dimensional diagram state is shown in Figure 4. According to the particularity of the position, the symmetry can simplify the three-dimensional model paraboloid to be formed by two-dimensional parabola rotation. In the model, a three-dimensional three-dimensional space rectangular coordinate system can be established with the reference spherical center C as the origin. In the simplified two-dimensional model, a plane rectangular coordinate system can be established with the reference spherical center C as the origin, assuming that the ideal two-dimensional parabolic equation is $y = a \cdot x^2 + c$.

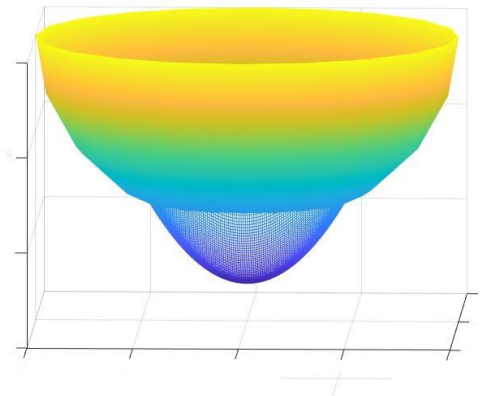

Figure4 3D schematic diagram of ideal paraboloid and reference sphere under problem 1

From the condition that the diameter of the approximate rotating paraboloid (working paraboloid) is 300 meters, the loss function can be calculated according to the sum of the squares of the distances between the reference sphere and the ideal paraboloid. The difference is a loss function. The Monte Carlo method is used to simulate the optimal ideal parabola iteratively, and then convert it into a three-dimensional ideal parabola to check whether the distance between adjacent nodes meets the constraints.

## 4.2 Analysis of Problem 2

The general idea of problem 2 is the extension of the ideal situation in the special state of problem 1, that is, when the observed celestial body S $\alpha = 36.795°$, $\beta = 78.169°$ is located in a non-trivial position, it is necessary to determine the ideal paraboloid, input the coordinates of the main cable node in the given conditions into the computer, and use the computer to program Python. The language can draw a schematic diagram of the point, see Figure 5.

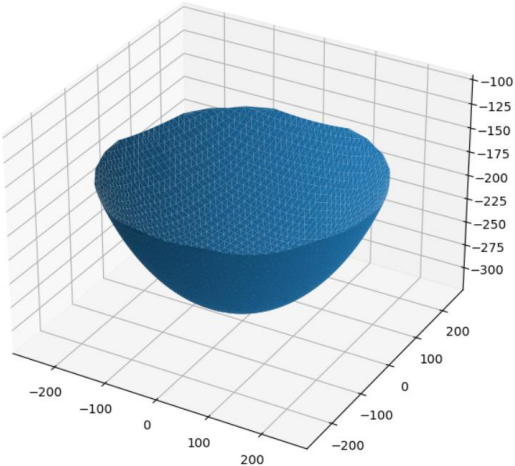

Figure 5 3D schematic diagram of discrete points of main cable node coordinates

The problem can be solved by means of the model established in Problem 1, by means of coordinate transformation after establishing the space rectangular coordinate system. After the transformation, the ideal paraboloid can be obtained under the condition that the position of the celestial body is limited by the computer simulation method, and the adjustment model of the reflective panel can be established. The working area is the red area, see the schematic diagram in Figure 6.

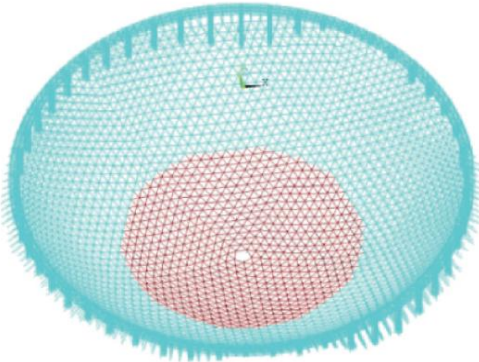

Figure 6 3D schematic diagram of the position comparison between the working paraboloid (red area) and the reference sphere (blue area)

## 4.3 Analysis of Problem 3

The third question is mainly based on the second question, and requires us to establish an evaluation index for the working paraboloid and the original spherical surface required in the second question to reflect the acceptance ratio of the feed cabin.

According to the data given in Annex 3, the fact that the area of each mirror is different is obtained through analysis, but the receiving range of the effective area of the feed cabin is only 1 meter, so it can be considered that as long as the signal reflected by one reflector can If it is received by the feed cabin, the effect is better.

Since the observed celestial body S is quite far away from the Chinese FAST radio telescope on the earth, the signal emitted by the observed celestial body can be approximately regarded as a parallel signal, and it is assumed that the signal is specularly reflected on the reflector. According to the second problem, the overall The coordinate axis transformation can transform the straight line SC axis to the Z axis of the coordinate system, so that the amount of calculation will be greatly reduced, and then judging whether the signal reflected by the reflector can be received by the feed cabin, the adjusted feed can be obtained. The source cabin acceptance ratio is compared with the acceptance ratio of the reference reflective sphere, and a conclusion is drawn.

# 5 Model establishment and solution

## 5.1 Establishment and solution of the problem-model

### 5.1.1 Model establishment

The celestial body to be measured S is located directly above the reference sphere. In this case $\alpha = 0°$, $\beta = 90°$, it is an ideal state. At this time, the three points of the observed celestial body S, the sphere center C of the reference sphere, and the feed cabin P are on the same straight line. The particularity here is that when the observed celestial body S is directly above the entire three-dimensional model is symmetrical about the vertical line SC, so according to the particularity of the position, the three-dimensional model paraboloid can be simplified to a two-dimensional parabolic rotation by symmetry Formed, a simplified diagram is shown in Figure 7.

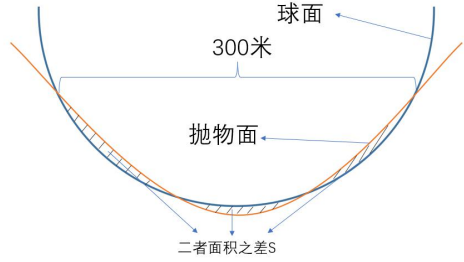

Fig.7 Simplified two-dimensional schematic diagram of ideal paraboloid and reference sphere under problem 1

Considering that in practical engineering applications, the changes of the main cable node and the entire reflector system need to be kept at a low level, so in problem 1, the purpose should be to ensure that the focus of the ideal paraboloid is the feed cabin P as much as possible. to reduce all node changes. According to the requirements, here is the situation where the ideal paraboloid needs to be solved, so the moving distance of the relevant points can be converted into the distance between the corresponding surfaces. In the two-dimensional case of the simplified model, the solution is based on the area S enclosed by the two lines. Make the area S as small as possible.

In the model, a three-dimensional solid space Cartesian coordinate system can be established with the reference spherical center C as the origin. In the simplified two-dimensional model, a plane rectangular coordinate system can be established with the reference sphere center C as the origin, assuming that the ideal two-dimensional parabola equation is $y = a \cdot x^2 + c$, and the reference sphere is simplified into a reference sphere semicircle, where according to the given condition R About 300.4 meters, the radius of the focal plane is the difference between the radius of the two concentric spheres $x^2 + y^2 = R^2$ of the reference sphere $F = 0.466R$. From the definition of the parabola and related properties, the following four formulas can be obtained to establish the model:

$$\begin{cases} \Gamma_1 : y_1 = ax_1^2 + c \\ \Gamma_2 : x_2^2 + y_2^2 = R^2 \\ F = 0.466R \\ R - F = \frac{1}{4a} + c \end{cases} \quad (1)$$

Then proceed to the solution process of the model.

### 5.1.2 The solution of the model

From the condition that the diameter of the approximate rotating paraboloid (working paraboloid) is 300 meters, the loss function can be calculated according to the sum of the squares of the distances between the reference sphere and the ideal paraboloid. The difference is the loss function, which is a representation of an integral.

The optimal ideal parabola can be solved step by step by approximating the micro-element method by the step path of computer simulation. Here, the method of definite integration can be used originally, but according to the situation of discrete points, the method of adding microelements is considered for simulation. Furthermore, since the distance between two straight lines is an absolute value function, it is not convenient to calculate the derivation, etc. operation, hereby changed to solve for the square of the difference between the curve distances. After solving the two-dimensional ideal parabola, convert it into a three-dimensional ideal parabola to check whether the distance between adjacent nodes satisfies the constraints. The following are the optimization conditions:

$$\min S = \int_{-150}^{150} (y_1 - y_2)^2 \, dx$$
$$= \int_{-150}^{150} \left( ax^2 + c - \sqrt{R^2 - x^2} \right)^2 dx \qquad (2)$$
$$= \lim_{\lambda \to 0} f(\xi_i) \Delta x_i$$

By writing a computer program code to solve it, the ideal parabola coefficient a is about 0.0017809, and c is about -300.79084, and the rotation is converted into a three-dimensional ideal paraboloid model, and the final ideal paraboloid equation is:

$$0.0017809(x^2 + y^2) = z + 300.79084 \qquad (3)$$

According to the relevant main cable node data, the obtained ideal paraboloid is compared with the main cable node coordinates in the original reference spherical surface. The requirement for the expansion and contraction range to be between plus and minus 0.6 meters, that is, the obtained ideal paraboloid effect is better.

## 5.2 The establishment and solution of the second problem model

The second problem is an extension of the ideal situation in the special state in the first problem, that is, when the observed celestial body S is located in a non-trivial position, the ideal paraboloid needs to be determined, and the specific position is α = 36.795°, β = 78.169°. In this case, the feed cabin can move on the focal plane, but the observed celestial body S, the reference sphere center C and the feed cabin P are still in a straight line, and the illumination area will be offset according to the position of the feed cabin , is no longer directly below the reference sphere under the first condition, but deforms into a paraboloid at the corresponding offset position to obtain a better reflection effect. The relevant schematic diagram is shown in Figure 8.

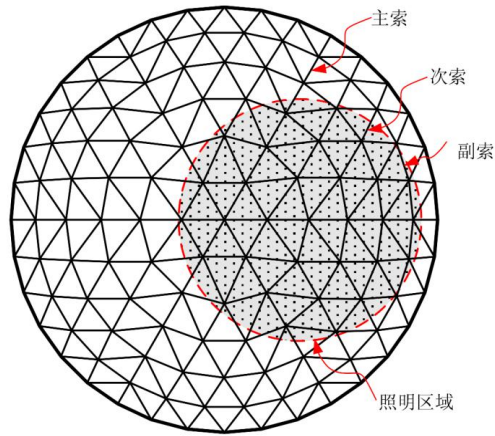

Fig. 8 Two-dimensional bird's-eye view of the ideal parabolic area and the reference sphere[3]

The main idea of solving the second problem is to select each discrete node position and each reflective panel as the variation between each main cable node, reflective panel and reflective surface, and ideal paraboloid, while the ideal paraboloid and the adjusted reflective surface need to be After the calculation is obtained, it is used as a reference to adjust each actuator to expand and contract to meet the goal of making the reflecting surface as close to the ideal paraboloid as possible.

Fig. 9 Schematic diagram of the geometric shape of the cable net under the working condition

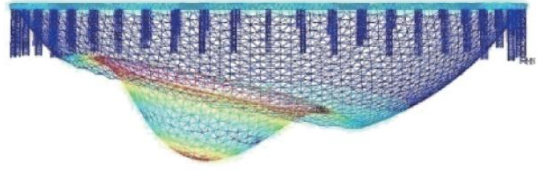

of the working paraboloid (exaggerated deformation in some areas) [2]

Note that in problem 2, the only change relative to problem 1 is the orientation of the celestial body. The center C of the reference sphere can be used as the origin, and the line connecting the observed celestial body S and the center of the reference sphere C is the Z axis to establish a space Cartesian coordinate system. In such a situation, the model and method obtained in the first question can be used to map all the main cable nodes and the relevant coordinates of the actuators in the data according to the transformation rules of the three-dimensional space to carry out new coordinate mapping.

According to the coordinate transformation model, the relevant derivation process is as follows:

1. First rotate the reference point around the Z axis by 0°, the rotation matrix is:

$$Rz = \begin{bmatrix} 1 & 0 & 0 \\ 0 & 1 & 0 \\ 0 & 0 & 1 \end{bmatrix} \quad (4)$$

2. Rotate the θ angle counterclockwise around the Y axis, and the CS direction vector is. Find the offset Cy' axis, so that y'⊥CS, CS' is the offset axis CX' of C, and move in the ZOX plane, then set the CX' The direction vector is . Then there are

$$\cos\beta\cos\alpha + a\sin\beta = 0$$
$$a = -\frac{\cos\beta\cos\alpha}{\sin\beta} \quad (5)$$

$$1 - b\frac{\cos\beta\cos\alpha}{\sin\beta} = 0$$
$$b = \frac{\sin\beta}{\cos\beta\cos\alpha} \quad (6)$$

CS'⊥CX', CX'⊥CS, then the plane CSS' is perpendicular to CX'.

Let the angle between CS' and CS be γ, then we have

$$\cos\gamma = \frac{\cos\beta\cos\alpha + \frac{\sin^2\beta}{\cos\beta\cos\alpha}}{\sqrt{1 + \frac{\sin^2\beta}{\cos^2\beta\cos^2\alpha}}} = \frac{\cos^2\beta\cos^2\alpha + \sin^2\beta}{\sqrt{\sin^2\beta + \cos^2\beta\cos^2\alpha}} \tag{7}$$

Let the angle between CS' and CX be θ', then we have

$$\cos\theta' = \frac{\cos\beta\cos\alpha}{\sqrt{\sin^2\beta + \cos^2\beta\cos^2\alpha}} \tag{8}$$

$$\theta = \frac{\pi}{2} - \theta' \tag{9}$$

Because $\theta = \frac{\pi}{2} - \theta'$, there is

$$R_y = \begin{bmatrix} \cos\left(\frac{\pi}{2} - \theta'\right) & 0 & -\sin\left(\frac{\pi}{2} - \theta'\right) \\ 0 & 1 & 0 \\ \sin\left(\frac{\pi}{2} - \theta'\right) & 0 & \cos\left(\frac{\pi}{2} - \theta'\right) \end{bmatrix}$$

$$= \begin{bmatrix} \sin\theta' & 0 & -\cos\theta' \\ 0 & 1 & 0 \\ \cos\theta' & 0 & \sin\theta' \end{bmatrix} \tag{10}$$

3. Rotate 2π-r counterclockwise around the X axis, then we have

$$R_x = \begin{bmatrix} 1 & 0 & 0 \\ 0 & \cos(2\pi - \gamma) & \sin(2\pi - \gamma) \\ 0 & -\sin(2\pi - \gamma) & \cos(2\pi - \gamma) \end{bmatrix}$$

$$= \begin{bmatrix} 1 & 0 & 0 \\ 0 & \cos\gamma & -\sin\gamma \\ 0 & \sin\gamma & \cos\gamma \end{bmatrix} \tag{11}$$

The final formula for R is

$$R_x R_y R_z = \begin{bmatrix} \sin\theta' & 0 & -\cos\theta' \\ -\sin\gamma\cos\theta' & \cos\gamma & -\sin\gamma\sin\theta' \\ \cos\gamma\cos\theta' & \sin\gamma & \cos\gamma\sin\theta' \end{bmatrix} \tag{12}$$

According to the above coordinate transformation formula, we import all relevant data into the computer, use MATLAB and Python programming, and calculate and solve the data such as the coordinates of new relevant points. The specific value of the transformation matrix obtained by solving is:

$$R = \begin{bmatrix} 0.9862209828 & 0 & -0.1654332887 \\ -0.02031535173 & 0.9924313308 & -0.1211087944 \\ 0.1641811789 & 0.1228008697 & 0.9787566025 \end{bmatrix} \quad (13)$$

$$R^\top = \begin{bmatrix} 0.9862209828 & -0.02031535173 & 0.1641811789 \\ -5.421010862e-20 & 0.9924313308 & 0.1228008697 \\ -0.1654332887 & -0.1211087944 & 0.9787566025 \end{bmatrix} \quad (14)$$

According to the solution process of the model in Problem 1, under the transformed coordinates in the obtained new space rectangular coordinate system, the equation of the ideal paraboloid under the condition of the specific celestial body position in Problem 2 can be simplified to the equation of the ideal paraboloid as:

$$\Gamma 3 : \begin{cases} x = 0.9862 x_1 - 0.1654 z_1 \\ y = -0.0203 x_1 + 0.9924 y_1 - 0.1211 z_1 \\ z = 0.1641 x_1 + 0.1228 y_1 + 0.9787 z_1 \\ 0.0017809 \left( x^2 + y^2 \right) = 2 + 300.79084 \end{cases} \quad (15)$$

Its vertex coordinates are(-49.3842 -36.9373 -294.4010).

According to the ideal paraboloid equation, we next fit the relevant coordinates when the actuator is not stretched or retracted in the reference state. We found a problem in the process of related experimental data processing: there are three coordinate points of some main cable nodes, the lower end point of the actuator (ground anchor point), and the upper end point (top) of the reference state, which cannot be used effectively. Fitted into a line, and further, some of the data gaps are large.

For example, the main index node with node number D69 has the coordinates of

$$M_{D69} = (-20.008, -109.112, -279.167)$$

The coordinates of the lower end point of the actuator (the ground anchor point) are

$$D_{D69} = (-20.809, -113.48, -286.222)$$

And the coordinates of the upper endpoint (top) in the reference state are

$$U_{D69} = (-20.676, -112.752, -284.385)$$

It can be verified that the three points are not collinear.

$$\begin{cases} \boldsymbol{v_{DU}}|_{x=1} = (1, 5.474, 13.812) \\ \boldsymbol{v_{UM}}|_{x=1} = (1, 5.449, 7.811) \end{cases}$$

It can be clearly seen that the components of the two direction vectors on the plane xOy are relatively consistent, but their components on the plane xOz are quite different.

In view of this, we choose the direction of the line connecting the actuator and the center of the reference sphere as the reference direction, and adjust the movement of the main cable node and the actuator to make it move to the ideal paraboloid as much as possible. According to the computer program design and operation model, using the formula:

$$\begin{cases} \boldsymbol{v} = \dfrac{\overrightarrow{UD}}{|UD|} \\ \dfrac{x - M_x}{\boldsymbol{v}_x} = \dfrac{y - M_y}{\boldsymbol{v}_y} = \dfrac{z - M_z}{\boldsymbol{v}_z} \\ z = a(x^2 + y^2) + c \end{cases} \quad (16)$$

Checking the relevant information again, according to the constraints, it is found that the solution formula for the change of the adjusted point relative to the point before adjustment is as follows:

$$dist = \sqrt{(N_x - M_x)^2 + (N_y - M_y)^2 + (N_z - M_z)^2} \quad (17)$$

We calculated the distance before and after adjustment for a total of 689 points within a diameter of 300 meters, and found that the minimum value of the before and after change distance was -0.3993m (node D27), and the highest value was 0.2833m (node B12).

The formula for solving the radial expansion and contraction of the top end of the range actuator is:

$$\begin{cases} dist(P_{id1}, P_{id2}) = \sqrt{\sum_{i}^{x,y,z}(P_{id1,i} - P_{id2,i})^2} \\ ratio(id1, id2) = \dfrac{|dist(N_{id1}, N_{id2}) - dist(M_{id1}, M_{id2})|}{dist(M_{id1}, M_{id2})} \end{cases} \quad (18)$$

Among all 3864 edges that satisfy the diameter of 300 meters, 2484 edges do not satisfy the change ratio of less than 0.07%, and the edge with the largest change is 0.62%.

In this case, the desired ideal paraboloid cannot meet the actual engineering requirements, and our point needs to be improved to make it closer to the ideal paraboloid while keeping the constraints.

The iterative adjustment method based on the least squares optimization model can be used, and the distance is used as the error to find the sum of the squares of the residuals of each node after adjusting the position, and perform the calculation. 0.6, and the variation of the distance between adjacent nodes is less than 0.07%. We use computer simulation to find the case where the residual sum of squares is the smallest under this constraint, and obtain the solution that is closest to the ideal paraboloid.

The constraint condition is mainly the change range of the distance between adjacent nodes, which requires our working paraboloid to approach the paraboloid as much as possible while satisfying the small change range of the distance between adjacent nodes. The range is relatively small, so the problem should not be considered from a single point, but should be considered from the whole surface. The solution is to make all the main cable nodes on the original reference sphere approximate the ideal paraboloid proportionally, and keep approaching until reaching Up to the upper bound of the allowable length variation.

$$\begin{cases} proportion = \dfrac{\overrightarrow{R_{id}M_{id}}}{\overrightarrow{N_{id}M_{id}}} \\ \max\left(ratio\left(id1,id2\right)\right) <= 0.07\% \\ \max proportion \end{cases} \tag{19}$$

Through computer simulation, we can obtain that the optimal ratio that can be achieved is about 0.89 under the constraint of guaranteeing a ratio of 0.07%.

Table 2 The relationship between the acceptance ratio of the feed cabin and the maximum variation range

| $proportion$ | … | 0.87 | 0.88 | 0.89 | 0.9 | … |
|---|---|---|---|---|---|---|
| $\max(ratio)$ | … | 0.077% | 0.071% | 0.065% | 0.060% | … |

After optimization, the working paraboloid under this problem has been determined, and then the coordinate points are mapped to this working paraboloid according to the relevant coordinate transformation steps above.

The change range of the point after re-adjustment relative to the point before adjustment is [-0.0443m, 0.0315m], and the maximum radial expansion and contraction of the actuator tip is 0.065%.

Finally, we store the results of the main cable node number and position coordinates obtained by the transformation and the results of the expansion and contraction of each actuator into the "result.xlsx" file.

Table 5 The number and position coordinates of the main cable node obtained by solving problem 2 after transformation

| Node Number | X Coordinate (m) | Y Coordinate (m) | Z Coordinate (m) |
|---|---|---|---|
| A0 | -0.00027625 | -0.000206619 | -300.4016468 |
| B1 | 6.108932021 | 8.407743517 | -300.2124881 |
| C1 | 9.883301288 | -3.210389947 | -300.2157062 |
| D1 | -0.001164183 | -10.39211915 | -300.2271401 |
| E1 | -9.884870582 | -3.212474613 | -300.230942 |
| A1 | -6.107932981 | 8.406931642 | -300.2209041 |
| A3 | 0.001492201 | 16.81860636 | -299.9200043 |

| …… | …… | …… | …… |
| D267 | -130.0763112 | -147.6225716 | -227.0389699 |
| D268 | -139.58752 | -139.8671354 | -226.3004543 |
| D269 | -148.8352568 | -131.8880491 | -225.1997059 |

## 5.3 Establishment and solution of the three-model problem

The third question is mainly based on the second question, and requires us to establish an evaluation index for the working paraboloid and the original spherical surface required in the second question to reflect the acceptance ratio of the feed cabin.

According to the data given in question A of the competition, among the 4300 mirrors, the minimum mirror area is 51.358m$^2$, and the main cable node number is: A0-C1-D1; while the maximum mirror area is 66.772m$^2$, the main cable node number is: A221 -A222-A243, the average size of all mirrors is 55.919m$^2$. The receiving range of the feed cabin is only 1 meter, so the size of the mirror is much larger than that of the feed cabin. In such a case, it can be considered that there is no need to consider the acceptance ratio of the reflected signal of a single reflector, as long as the signal reflected by one reflector can be received by the feed cabin, then a better effect can be achieved.

Since the observed celestial body S is quite far away from the FAST radio telescope on the earth, the signal emitted by the observed celestial body S can be approximately regarded as a parallel signal, and it is assumed that the signal is specularly reflected on the reflector. , the overall coordinate axis transformation is carried out, so the straight line SC can be transformed to the Z axis, which can effectively reduce the calculation amount and complexity of the model.

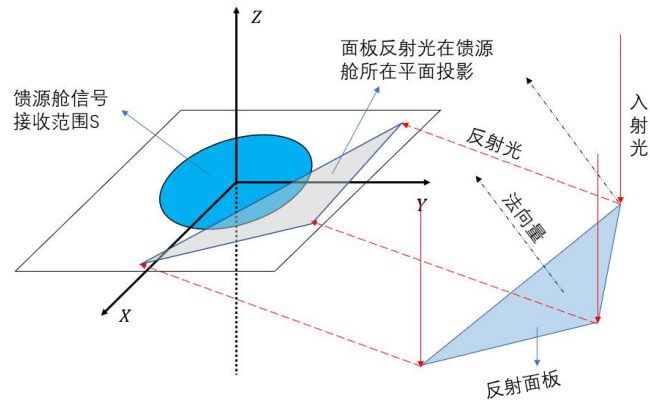

Figure 10 Simple schematic diagram of specular reflection model

According to the law of reflection [4]: When light strikes an interface, the incident ray and the reflected ray make the same angle. When light is incident on the interface of different media, reflection and refraction occur. The reflected ray is in the same plane as the incident ray and the normal. The reflected line and the incident line are separated from the normal line, and the angle between them and the interface normal line (called the angle of incidence and the angle of reflection, respectively) is equal, and the angle of reflection is equal to the angle of incidence.

It can be obtained that the angle between the reflected signal vector $\boldsymbol{n}$ and the plane normal vector $\boldsymbol{n}_0$ is complementary to the angle between the incident signal vector $\boldsymbol{n}_z$ and the plane normal vector $\boldsymbol{n}_0$, so that the expression can be listed

$$\cos <\boldsymbol{n}, \boldsymbol{n}_0> = \cos <\boldsymbol{n}_0, \boldsymbol{n}_z> \tag{20}$$

Knowing that these three vectors are on the same plane, that is to say, the mixed product of these three vectors is 0, we can get

$$\boldsymbol{n} \cdot (\boldsymbol{n}_0 \times \boldsymbol{n}_z) = 0 \tag{20}$$

Let us assume that all three are unit vectors.

$$\begin{cases} \boldsymbol{n}_z = (0,0,1) \\ \boldsymbol{n}_0 = (x_0, y_0, z_0) \\ \boldsymbol{n} = (x, y, z) \\ |\boldsymbol{n}| = |\boldsymbol{n}_0| = 1 \end{cases} \tag{21}$$

From this, the following equation can be obtained

$$\begin{cases} x_0 x + y_0 y + z_0 z = 1 \\ \det\begin{pmatrix} 0 & 0 & 1 \\ x_0 & y_0 & z_0 \\ x & y & z \end{pmatrix} = 0 \\ x^2 + y^2 + z^2 = 1 \end{cases} \tag{22}$$

Solutions have to

$$\begin{cases} x = \dfrac{2 x_0 z_0}{x_0^2 + y_0^2 + z_0^2} \\ y = \dfrac{2 y_0 z_0}{x_0^2 + y_0^2 + z_0^2} \\ z = \dfrac{z_0^2 - x_0^2 - y_0^2}{x_0^2 + y_0^2 + z_0^2} \end{cases} \tag{23}$$

The next step is how to judge whether the signal reflected by the reflector can be received by the feed cabin.

According to the requirements of the subject, the range of the feed cabin is only 1 meter. According to the requirements of the subject, the center of the receiving plane of the feed cabin can only move on a spherical surface (focal plane) that is concentric with the reference spherical surface. It can be concluded that the feed cabin can receive The equation for the resulting signal is:

$$\begin{cases} R = 300.4m \\ F = 0.466R \\ z = -(R-F) \\ x^2 + y^2 <= 1 \end{cases} \qquad (24)$$

Therefore, it is only necessary to consider the distribution of the reflected signal on the plane xOy where the feed cabin is located.

Since it is assumed that the incident signal is parallel and the angle formed with the reflector is always the same, it can be concluded that, according to the linear invariance of space, the graph of the reflected signal on the plane where the feed cabin is located is also a triangle. From the plane reflection vector obtained above, we can derive some of the following formulas.

Suppose the three points on the reflective panel are $A, B, C$, and the three vertices of the triangle formed by the projection of the signal reflected by the reflective panel on the plane of the feed bin are $A', B', C'$, then

$$\begin{cases} \boldsymbol{v}_1 = \overrightarrow{AB}, \boldsymbol{v}_2 = \overrightarrow{AC} \\ \boldsymbol{n}_0 = \boldsymbol{v}_1 \times \boldsymbol{v}_2 \\ \boldsymbol{n} = f(\boldsymbol{n}_0) = (n_x, n_y, n_z) \\ P_{x,y,z}' = \begin{cases} \dfrac{x - P_x}{n_x} = \dfrac{y - P_y}{n_y} = \dfrac{z - P_z}{n_z} \quad (P = A, B, C) \\ z = -(R-F) \end{cases} \end{cases} \qquad (25)$$

Wherein f in formula (25-3) is the calculation result in formula (21).

Then calculate whether the triangle and the circle have a common area.

If there is a common area between a circle and a triangle, there are only the following two situations: the center of the circle is inside the triangle, and the circle and the triangle have an intersection, as shown in Figure 11:

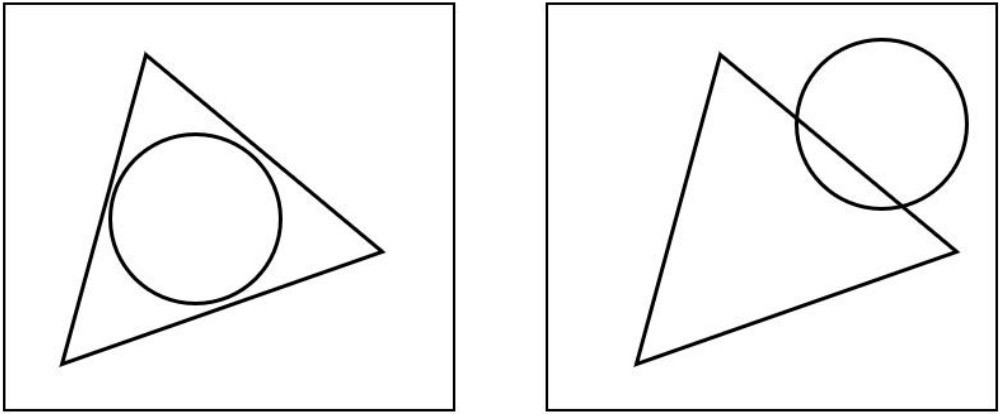

Figure 11 A simple schematic diagram of a circle inside and outside a triangle

①The circle is inside the triangle:

The Triangle Interior Algorithm is used here, that is, the determination theorem of the interior points of a triangle. The following is a brief introduction to the theorem:

To judge whether a given point $P$ is in $\triangle ABC$ or not, you can do the following, let

$$\begin{cases} \boldsymbol{v} = \overrightarrow{OP} \\ \boldsymbol{v}_0 = \overrightarrow{OA} \\ \boldsymbol{v}_1 = \overrightarrow{AB} \\ \boldsymbol{v}_2 = \overrightarrow{AC} \end{cases} \quad (26)$$

Then according to the vector decomposition theorem, we have

$$\boldsymbol{v} = \boldsymbol{v}_0 + a\boldsymbol{v}_1 + b\boldsymbol{v}_2 \quad (27)$$

in

$$\begin{cases} a = \dfrac{\det(\boldsymbol{v}\boldsymbol{v}_2) - \det(\boldsymbol{v}_0\boldsymbol{v}_2)}{\det(\boldsymbol{v}_1\boldsymbol{v}_2)} \\ b = -\dfrac{\det(\boldsymbol{v}\boldsymbol{v}_1) - \det(\boldsymbol{v}_0\boldsymbol{v}_1)}{\det(\boldsymbol{v}_1\boldsymbol{v}_2)} \end{cases} \quad (28)$$

if satisfied

$$\begin{cases} a, b > 0 \\ a + b < 1 \end{cases} \quad (29)$$

then the point $P$ is inside $\triangle ABC$.

Here you only need to determine whether the center of the circle is inside the triangle.

②Circles and triangles have a common point:

There is a common point, that is, the intersection of the circle and the line segment. We can calculate the line segment equation by substituting it into the circle.

$$\begin{cases} x_1, x_2 = \begin{cases} y = kx + b \, (defined \; by \; line) \\ x^2 + y^2 = 1 \end{cases} \\ \exists \, x, \; s.t. \; x_i \in dom(X_{line}) \end{cases} \quad (30)$$

Conclusion can be made:

The acceptance ratio of the reflector signal is defined as the ratio of the number of surfaces that the signal can reach the feed cabin to the number of all surfaces, that is, the acceptance ratio *efficiency* is defined as follows:

$$efficiency = 100\% \times \frac{\forall \triangle, S(\odot \cap \triangle) > 0}{tot \triangle} \qquad (31)$$

In the original positive sphere, there are a total of 4300 reflectors, of which 832 can reflect the signal, accounting for the working parabola established based on the data of our problem 2, there are a total of 1288, of which 805 can reflect the signal, It accounts for 62.5%. Compared with the two, it has increased by 224%, and the effect is better.

# 6 MODEL EVALUATION

## 6.1 Advantages of the model

In the first problem, the problem in the three-dimensional space is simplified to a two-dimensional problem, which is convenient for solving and calculation. When considering the constraints on points in the problem, the point is replaced by a surface, which greatly simplifies the problem. Restrictions on points serve as validation conditions rather than constraints, allowing the number of variables in the problem to be reduced.

In problem 2, when solving the model, we keenly notice the correlation between problem 2 and problem 1, and use the space transformation matrix to transform the point solved in problem 2 to the position of the point in problem 1, which is a good implementation. The reuse of calculation code and parabolic formula makes the model have good generalization. When analyzing the constraints, the target signal is not restricted to a single point, but the transformation of the entire surface is considered, which reduces the complexity.

In problem 3, when solving the reflection problem, the simulation of the signal line is not considered, but the real reflection surface is solved through the derivation of the mathematical formula, so that the order of magnitude of the simulation changes from the signal unit to the reflection surface, and attention is paid to the feed source. When calculating the receiving range of the cabin, the three-dimensional three-dimensional problem is reduced to a two-dimensional plane for solving, which reduces the amount of calculation and makes the results easy to verify.

## 6.2 Disadvantages of the model

In problem 1, only the two-dimensional situation is considered when solving the loss function, and the loss function that does not consider the three-dimensional and two-dimensional loss functions may not be linear. The numerical solution, not the analytical solution, is used in the solution.

In the second problem, when solving the working paraboloid, a simple sphere and an ideal paraboloid are used for interpolation. Later, we consider that the use of a piecewise paraboloid may have a better effect on the model.

In question 3, the amount of signal reflected by a single panel is not considered when calculating the reflectivity, and the data obtained by simply judging with a 0-1 variable may not be accurate.

### 6.3 Improvements to the model

When solving the working paraboloid, various practical situations should be considered, and different models should be established to evaluate its indicators, including but not limited to second-order fitting, linear interpolation, piecewise fitting, least squares and other methods. When establishing the loss function, the constraints of the problem should be fully considered, and the change of the cable length should be the main limiting factor, and the key should be grasped. There may be better results.

### 6.4 Generalization of the model

In this modeling, a large number of symmetry methods are used to achieve data dimensionality reduction. This idea can be applied to many practical engineering problems, and in general, the time complexity can be reduced by an order of magnitude.

It is not difficult to find that the problem is not limited by the position of the celestial body when solving this problem. If the azimuth and elevation angles of the celestial body change, it only needs to change a few parameters. Effect.

The method of solving the reflected signal in this question has a strong generalization. It uses the spatial invariance of the reflecting surface to transform the problem from a three-dimensional space to solving the intersection on a two-dimensional plane. According to relevant data, this algorithm is currently used in the signal tracking engine. Some have quite a few applications.